\documentstyle[12pt]{article}

\begin{document}
\title {A SIMPLE PROOF OF A THEOREM ON $(2n)$-WEAK AMENABILITY }
\author{M. S. Moslehian \\F. Negahban}
\date{}
\maketitle{}
{\bf Abstract.}{\small~ A simple proof of $(2n)$-weak amenability of the triangular Banach algebra ${\cal T}= \left [ \begin{array}{cc}{\cal A}&{\cal A}\\0&{\cal A} \end{array}\right ]$ is given where ${\cal A}$ is a unital $C^*$-algebra.}
\newpage
\begin{center}
{\bf 1. Introduction}
\end{center}

The topological cohomology groups provide us some significant information about Banach algebras such as their amenability, contractibility, stability, and singular extensions.[3]

Suppose that ${\cal A}$ and ${\cal B}$ are unital Banach algebras. ${\cal M}$ is a unital Banach ${\cal A}-{\cal B}-$module whenever it is simultaneously a Banach space, a left ${\cal A}$-module and a right ${\cal B}$-module satisfying $a(mb)=(am)b, ~1_{\cal A}m=m1_{\cal B}$ and $\parallel axb\parallel \leq \parallel a \parallel \parallel x \parallel \parallel b\parallel$. Then ${\cal T}= \left [ \begin{array}{cc}{\cal A}&{\cal M}\\0&{\cal B} \end{array}\right ] = \{ \left [ \begin{array}{cc}a&m\\0&b \end{array}\right ]; ~a\in {\cal A}, m\in {\cal M}, b\in {\cal B} \}$ with the usual $2\times 2$ matrix addition and formal multiplication equipped with the norm $\parallel \left [ \begin{array}{cc}a&m\\0&b \end{array}\right ]\parallel = \parallel a \parallel + \parallel m \parallel + \parallel b \parallel$ is said to be a triangular Banach algebra.

Note that the dual ${\cal M^*}$ of ${\cal M}$ together with the actions $(\phi a)(x)=\phi(ax)$ and $(b\phi)(x)=\phi(xb)$ is a Banach ${\cal B}-{\cal A}-$module. Similarly the $(2n)$-th dual ${\cal M}^{(2n)}$ of ${\cal M}$ is a Banach ${\cal A}-{\cal B}-$module and the $(2n-1)$-th dual ${\cal M}^{(2n-1)}$ of ${\cal M}$ is a Banach ${\cal B}-{\cal A}-$module. In particular, ${\cal A}^{(n)}$ is a Banach {\cal A}-bimodule.

The notion of $n$-weak amenability was introduced by Dales, Ghahramani and Gronb$\ae$ck [1]. If the first topological cohomology group $H^1({\cal A},{\cal A}^{(n)})$, i.e the quotient $\{\delta: {\cal A} \longrightarrow {\cal A}^{(n)}; \delta\\
{\rm ~ is~linear~ and ~}\delta(ab)=a\delta(b)-\delta(a)b\}$ modulo $\{\delta_x: {\cal A} \longrightarrow {\cal A}^{(n)}; \delta_x(a)=ax-xa, a\in {\cal A}, x\in {\cal A}^{(n)}\}$ is zero, then ${\cal A}$ is called $n$-weakly amenable. If for all $n, {\cal A}$ is $n$-weakly amenable, ${\cal A}$ is said to be permanently weakly amenable. For example, every $C^*$-algebra is permanently weakly amenable [1, Theorem 3.1].

Forrest and Marcoux investigated a relation between $n$-weak amenability of triangular Banach algebra ${\cal T}= \left [ \begin{array}{cc}{\cal A}&{\cal M}\\0&{\cal B} \end{array}\right ]$ and those of algebras ${\cal A}$ and ${\cal B}$. In particular, they proved permanently weak amenability of the triangular Banach algebra ${\cal T}= \left [ \begin{array}{cc}{\cal A}&{\cal A}\\0&{\cal A} \end{array}\right ]$ where ${\cal A}$ is a $C^*$-algebra.

In this paper we give a simple proof of the $(2n)$-weak amenability of the Banach algebra ${\cal T}= \left [ \begin{array}{cc}{\cal A}&{\cal A}\\0&{\cal A} \end{array}\right ]$ in which ${\cal A}$ is a $C^*$-algebra .[2, Proposition 4.5]

\begin{center}
{\bf 2.Main Result}
\end{center}
We start by reviewing some results on $(2n-1)$-weak amenability of triangular Banach algebras:

{\bf Theorem 2.1.} {\it ~Let ${\cal A}$ and ${\cal B}$ be unital Banach algebras and ${\cal M}$ be a unital Banach ${\cal A}-{\cal B}-$module. Let ${\cal T}= \left [ \begin{array}{cc}{\cal A}&{\cal M}\\0&{\cal B} \end{array}\right ]$ be the corresponding triangular Banach algebra. Let $n$ be a positive integer. Then $$H^1({\cal T},{\cal T}^{(2n-1)})\simeq H^1({\cal A},{\cal A}^{(2n-1)})\oplus H^1({\cal B},{\cal B}^{(2n-1)}).$$
It follows that ${\cal T}$ is $(2n-1)$-weakly amenable if and only if both ${\cal A}$ and ${\cal B}$ are.}[2, Theorem 3.7]

{\bf Corollary 2.2.} {\it ~${\cal T}= \left [ \begin{array}{cc}{\cal A}&{\cal A}\\0&{\cal A} \end{array}\right ]$
is $(2n-1)$-weakly amenable if ${\cal A}$ is a $C^*$-algebra.}

{\bf Proof.} Apply above theorem and  permanently weak amenability of ${\cal A}$ [1, Theorem 3.1].\\

{\bf Definition 2.4.} Let ${\cal M}$ be a Banach ${\cal A}-{\cal B}-$module and $\rho_{x,y}: {\cal M} \longrightarrow {\cal M}^{(2n)}$ be defined by $\rho_{x,y}(m)=xm-my$, where $x\in {\cal A}^{(2n)}, ~m\in {\cal M}$ and $y\in {\cal B}^{(2n)}$. Recall that if ${\cal A}$ is a Banach algebra, so is ${\cal A}^{(2n)}$ equipped with the (first) Arens product $\Gamma_1 \circ \Gamma_2=w^*-\lim_i\lim_ja_ia_j$ where $\{a_i\}$ and $\{a_j\}$ are nets in ${\cal A}^{(2n-2)}$ converging in weak$^*$-topology to $\Gamma_1, \Gamma_2\in {\cal A}^{(2n-2)}$, respectively. We could therefore define the centeralizer of ${\cal A}$ in ${\cal A}^{(2n)}$ as $Z_{\cal A}({\cal A}^{(2n)})=\{x\in {\cal A}^{(2n)}; xa=ax {\rm~for~ all~} a\in {\cal A}\}$ and the central Rosenblum operator on ${\cal M}$ with coefficients in ${\cal M}^{(2n)}$ as $ZR_{{\cal A},{\cal B}}({\cal M},{\cal M}^{(2n)})=\{\rho_{x,y}; x\in Z_{\cal A}({\cal A}^{(2n)}), y\in Z_{\cal B}({\cal A}^{(2n)})\}$. The later space is clearly a subspace of $Hom_{{\cal A},{\cal B}}({\cal M},{\cal M}^{(2n)})=\{\phi:{\cal M}\longrightarrow {\cal M}^{(2n)}; \phi(amb)=a\phi(m)b, {\rm ~for ~all~} a\in {\cal A}, m\in {\cal M}, b\in {\cal B}\}$.\\
The following theorem play a key role in the subject:

{\bf Theorem 2.5.} {\it Let ${\cal A}$ and ${\cal B}$ be unital Banach algebras and ${\cal M}$ be a unital Banach ${\cal A}-{\cal B}-$module. Let ${\cal T}= \left [ \begin{array}{cc}{\cal A}&{\cal M}\\0&{\cal B} \end{array}\right ]$ be the corresponding triangular Banach algebra. If $n$ is a positive integer and both ${\cal A}$ and ${\cal B}$ are $(2n)$-weakly amenable, then }$$H^1({\cal T},{\cal T}^{(2n)})=H^1({\cal M},{\cal M}^{(2n)})/ ZR_{{\cal A},{\cal B}}({\cal M},{\cal M}^{(2n)}).$$

We are ready to give our proof of Proposition 4.5 of [2]:

{\bf Theorem 2.6.} {\it ~${\cal T}= \left [ \begin{array}{cc}{\cal A}&{\cal A}\\0&{\cal A} \end{array}\right ]$
is $(2n)$-weakly amenable if ${\cal A}$ is a $C^*$-algebra.}

{\bf Proof.} For each $n$, the von Neumann algebra ${\cal A}^{(2n)}$ is unital with the unit denoted by $1_{{\cal A}^{(2n)}	}$. Clearly $ZR_{{\cal A},{\cal A}}({\cal A},{\cal A}^{(2n)}) \subseteq Hom_{{\cal A},{\cal A}}({\cal A},{\cal A}^{(2n)})$. For the converse, assume that $\phi\in Hom({\cal A},{\cal A}^{(2n)})$. Since $\phi(a)=\phi(a.1.1)=a\phi(a)$ and $\phi(a)=\phi(1.1.a)=\phi(1)a$, we have $\phi(1)\in Z_{\cal A}({\cal A}^{(2n)})$ and $1_{{\cal A}^{(2n)}}-\phi(1)\in Z_{\cal A}({\cal A}^{(2n)})$. Hence $\phi(a)=1_{{\cal A}^{(2n)}}.a-1_{{\cal A}^{(2n)}}.a+\phi(a)=\rho_{1_{{\cal A}^{(2n)}},1_{{\cal A}^{(2n)}}-\phi(a)}(a)$. 
It follows that $\phi=\rho_{1_{{\cal A}^{(2n)}},1_{{\cal A}^{(2n)}}-\phi(a)}\in ZR_{{\cal A},{\cal A}}({\cal A},{\cal A}^{(2n)})$. Therefore $ZR_{{\cal A},{\cal A}}({\cal A},{\cal A}^{(2n)}) =Hom_{{\cal A},{\cal A}}({\cal A},{\cal A}^{(2n)})$. Applying Theorem 2.5, we conclude that $$H^1({\cal T},{\cal T}^{(2n)})=H^1({\cal A},{\cal A}^{(2n)})/ ZR_{{\cal A},{\cal A}}({\cal A},{\cal A}^{(2n)})=0.$$

\end{document}